\documentclass[11pt]{article}
\usepackage{cite}
\usepackage{amsmath}
\usepackage{amsfonts}
\usepackage{graphicx}
\usepackage{subfigure}

\numberwithin{equation}{section}
\newcommand{\ep}{\varepsilon}

\newcommand{\p}{\partial}

\newcommand{\sech}{{\mathrm{sech}} }

\newcommand{\fl}[2]{\frac{#1}{#2}}
\newcommand{\be}{\begin{equation}}
\newcommand{\ee}{\end{equation}}
\newcommand{\nn}{\nonumber}

\newtheorem{theorem}{Theorem}[section]

\newtheorem{corollary}[theorem]{Corollary}
\newtheorem{lemma}[theorem]{Lemma}
\newtheorem{remark}[theorem]{Remark}

\begin{document}

\title{Two exponential-type integrators for the ``good" Boussinesq equation}

\author{Alexander Ostermann\thanks{Department of Mathematics,
University of Innsbruck, 6020 Innsbruck, Austria ({\tt alexander.ostermann@uibk.ac.at})}\,\,\,\,and
Chunmei Su\thanks{Department of Mathematics,
University of Innsbruck, 6020 Innsbruck, Austria ({\tt sucm13@163.com})}}
\date{}
\maketitle

\begin{abstract}
We introduce two exponential-type integrators for the ``good" Bousinessq equation.
They are of orders one and two, respectively, and they require lower regularity of the solution compared to the classical exponential integrators. More precisely, we will prove first-order convergence in $H^r$ for solutions in $H^{r+1}$ with $r>1/2$ for the derived first-order scheme. For the second integrator, we prove second-order convergence in $H^r$ for solutions in $H^{r+3}$ with $r>1/2$ and convergence in $L^2$ for solutions in $H^3$. Numerical results are reported to illustrate the established error estimates. The experiments clearly demonstrate that the new exponential-type integrators are favorable over classical exponential integrators for initial data with low regularity.
\end{abstract}


\pagestyle{myheadings}\thispagestyle{plain}

\section{Introduction}
Consider the ``good" Boussinesq (GB) equation \cite{boussinesq1872}
\be\label{GB}
z_{tt}+z_{xxxx}-z_{xx}-(z^2)_{xx}=0,
\ee
which was originally introduced as a model for one-dimensional weakly nonlinear dispersive waves in shallow water. Similar to the well-known Korteweg--de Vries (KdV) equation and the cubic Schr\"odinger equation, the GB equation is one of the important models describing the interaction between nonlinearity and dispersion. The GB equation has been widely applied in many areas, e.g., plasma, coastal engineering, hydraulics studies, elastic crystals and so on.

The GB equation and its various extensions have been extensively analyzed in the literature. For the well-posedness, we refer to \cite{fang1996, farah2009, farah2010, wang2013, oh2013improved, kishimoto2013sharp} and references therein. For the interaction of solitary waves, we refer to \cite{manoranjan1984, manoranjan1988}. Many numerical methods have been developed for solving the GB equation, such as finite difference methods (FDM) \cite{ortega1990, bratsos2007}, Petrov-Galerkin methods \cite{manoranjan1984}, mesh free methods \cite{dehghan2012}, Fourier spectral methods \cite{cheng2015, yan2016, de1991} and operator splitting methods \cite{zhang2017, zhang2018}. Regarding the numerical analysis for the
GB equation, nonlinear stability and convergence of some finite difference methods were studied in \cite{ortega1990}. Specifically, an explicit finite difference scheme was proved to converge quadratically in both space and time under the regularity assumption that $\fl{\partial^6 z}{\partial x^6}, \fl{\partial^4 z}{\partial t^4}$ are bounded and a severe time step restriction: $\Delta t=O(\Delta x^2)$ where $\Delta t$ and $\Delta x$ represent the discretization parameters in time and space, respectively. For a pseudospectral discretization with periodic boundary conditions, a second order temporal discretization was proposed and analyzed in \cite{de1991}, and full order convergence was proved in a weak energy norm: the $L^2$ norm in $z$ combined with the $H^{-2}$ norm in $z_t$ under a similar time step constraint $\Delta t=O(\Delta x^2)$. Due to the absence of a dissipation mechanism in the GB equation \eqref{GB}, it is more challenging to analyze the nonlinear error terms here than for a parabolic equation. The presence of a second-order spatial derivative in the nonlinear term brings an essential difficulty for numerical error estimates in a higher order Sobolev norm \cite{cheng2015}. The norm was strengthened in \cite{cheng2015}, where a second order temporal scheme was proposed and convergence was proved in a stronger energy norm: the $H^2$ norm in $z$ combined with the $L^2$ norm in $z_t$. Moreover, such a convergence is unconditional so that the severe time step condition $\Delta t=O(\Delta x^2)$ is avoided. However, it requires the solution to be smooth enough such that $z\in H^{m+4}$, $z_{tt}\in H^4$ and $\partial_t^4 z\in L^2$ to get an error of $O(\Delta t^2+\Delta x^m)$. Such constraints become very restrictive for computation, especially for rough solutions.

In this paper, we present two exponential-type Fourier integrators for the GB equation which enable us to lower the classical regularity assumptions and obtain a first-order and second-order convergence by requiring one and three additional derivatives, respectively. The exponential-type integrators are constructed based on the following strategy:
 \begin{enumerate}
 \item[1.] In the first step, we formulate the GB equation \eqref{GB} as a first-order equation in the complex domain via the transformation
     \[u=z-i\langle \p_x^2\rangle_c^{-1}z_t.\]

 \item[2.] In the second step, we rescale the first-order equation in time by considering the so-called \emph{twisted variable}
     \[w(t)=e^{it\p_x^2}u(t).\]
     This essential step will later enable us to treat the dominant term triggered by the nonlinearity in an exact way.

 \item[3.] Finally, we iterate Duhamel's formula in $w(t)$ and integrate the dominant interactions exactly.
 \end{enumerate}
The idea of twisting the variable is widely applied in the analysis of PDEs in low regularity spaces \cite{bourgain1993}. It was also widely applied in the context of numerical analysis for the Schr\"odinger equation \cite{ostermann2018low, knoller2018fourier}, the KdV equation \cite{hofmanova2017} and Klein-Gordon type equations \cite{baumstark2018, baumstark2018uniformly}.

For implementation issues, we impose periodic boundary conditions and
refer to \cite{farah2010, oh2013improved, kishimoto2013sharp} for the corresponding well-posedness results. For $m\in \mathbb{R}$, we define the Sobolev norm on $\Omega=(-\pi,\pi)$ by
\[\|f\|_m^2=\sum\limits_{k=-\infty}^\infty (1+k^2)^{m}
|\widehat{f}_k|^2,\quad \mathrm{where}\quad \widehat{f}_k=\fl{1}{2\pi}\int_\Omega f(x)e^{-ikx}dx.\]
Moreover, we denote by $H^m$ all the functions defined on $\Omega$ with finite norm $\|\cdot\|_m$. For $m=0$, the space is exactly $L^2$ and the corresponding norm is denoted as $\|\cdot\|$.

The rest of this paper is organized as follows. In Section 2, we introduce a scaling of the GB equation and give some preliminary notations and lemmas. The first- and second-order exponential-type integrators are constructed and analyzed in Sections 3 and 4, respectively. Numerical results are reported in Section 5, which illustrate the proved convergence results. Finally, some concluding remarks are drawn in Section 6.

\section{Scaling for low regularity exponential-type integrators}
Consider the GB equation
\be\label{GBp}
z_{tt}+z_{xxxx}-z_{xx}-(z^2)_{xx}=0,\quad x\in (-\pi, \pi), \quad t>0,
\ee
which can be reformulated as a first-order coupled system as
\begin{align*}
&\partial_t u=i\langle \p_x^2\rangle_c u-\fl{i}{4}\langle \p_x^2\rangle_c^{-1}\left[2c(u+\overline{v})+\p_{xx}(u+\overline{v})^2\right],\\
&\partial_t v=i\langle \p_x^2\rangle_c v-\fl{i}{4}\langle \p_x^2\rangle_c^{-1}\left[2c(\overline{u}+v)+\p_{xx}(\overline{u}+v)^2\right],
\end{align*}
where
\be\label{cdef}
\langle \p_x^2\rangle_c=\sqrt{\p_x^4-\p_x^2+c},\quad
u=z-i\langle \p_x^2\rangle_c^{-1}z_t,\quad
v=\overline{z}-i\langle \p_x^2\rangle_c^{-1}\overline{z_t},
\ee
and $c$ is a positive number. Since $z\in \mathbb{R}$ we immediately get $u=v$ and
\be\label{rec}
z=\fl{1}{2}(u+\overline{u}),\quad z_t=\fl{i}{2}\langle \p_x^2\rangle_c (u-\overline{u}).
\ee
Thus the GB equation \eqref{GBp} is equivalent to the following first-order equation in the complex domain
\be\label{NF}
\partial_t u=i\langle \p_x^2\rangle_c u-\fl{i}{4}\langle \p_x^2\rangle_c^{-1}\left[2c(u+\overline{u})+\p_{xx}(u+\overline{u})^2\right],\quad x\in(-\pi, \pi),\quad t>0,
\ee
with initial data
\[u(0)=z(0)-i\langle \p_x^2\rangle_c^{-1}z_t(0).\]
\begin{remark}
The linear term ``$cz$" was added and subtracted in equation \eqref{GBp} so that $u$ and $v$ in \eqref{cdef} are well-defined for functions with $\widehat{(z_t)}_0\neq 0$.
\end{remark}

We note that the leading term in $\langle \p_x^2\rangle_c$ is $-\p_x^2$. This motivates us to introduce the so-called \emph{twisted variable}
\[w(t)=e^{it\p_x^2}u(t).\]
In this new variable \eqref{NF} becomes
\be\label{wf}
\begin{split}
\partial_t w&=iAw-\fl{i}{4}e^{it\p_x^2}\langle \p_x^2\rangle_c^{-1}\left[2c(e^{-it\p_x^2}w+e^{it\p_x^2}\overline{w})+\p_x^2(e^{-it\p_x^2} w+e^{it\p_x^2}\overline{w})^2\right]\\
&=iAw-\fl{ic}{2}\langle \p_x^2\rangle_c^{-1}w-\fl{ic}{2}e^{2it\p_x^2}\langle \p_x^2\rangle_c^{-1}\overline{w}-\fl{i}{4}e^{it\p_x^2}\langle \p_x^2\rangle_c^{-1}\p_x^2
(e^{-it\p_x^2} w+e^{it\p_x^2}\overline{w})^2,
\end{split}
\ee
where $A=\p_x^2+\langle \p_x^2\rangle_c$. Thus by Duhamel's formula, we have
\be\label{Duhm}
\begin{split}
w(t_n+\tau)&=e^{i\tau A}w(t_n)-\fl{ic}{2}\langle \p_x^2\rangle_c^{-1}\int_0^\tau
e^{i(\tau-s)A}\left[w(t_n+s)+e^{2i(t_n+s)\p_x^2}\overline{w}(t_n+s)\right]ds\\
&\hspace{-10mm}-\fl{i\p_x^2}{4}\langle \p_x^2\rangle_c^{-1}\int_0^\tau e^{i(\tau-s)A}e^{i(t_n+s)\p_x^2}\left[e^{-i(t_n+s)\p_x^2}w(t_n+s)+e^{i(t_n+s)\p_x^2}
\overline{w}(t_n+s)\right]^2ds.
\end{split}
\ee

In the remainder of this paper we assume $r>1/2$ so that the well-known bilinear estimate holds:
\be\label{bi}
\|fg\|_r\le D_r\|f\|_r\|g\|_r,
\ee
where $D_r$ represents a positive constant depending on $r$. For simplicity of notation, we will employ the following definitions.

\emph{Notations:} Throughout the paper we will use the following notations
\[\psi_1(z)=\int_0^1 e^{zs}ds,\quad \psi_2(z)=\int_0^1s e^{zs}ds, \quad \mathrm{for} \quad z\in \mathbb{C}.\]
For $f(x)=\sum\limits_{k\in\mathbb{Z}}\widehat{f}_k e^{ikx}$, the regularization of $\p_x^{-1}$ and $\p_x^{-2}$ is defined through its action in Fourier space by
\[(\p_x^{-1})_k:=\left\{\begin{aligned}
&(ik)^{-1}\quad &\mathrm{if}\quad k\neq 0,\\
&0\quad &\mathrm{if}\quad k=0,
\end{aligned}\right.
\quad \mathrm{i.e.},\quad \p_x^{-1} f(x)=\sum\limits_{k\neq 0} \fl{\widehat{f}_k}{ik} e^{ikx},\quad
\p_x^{-2} f(x)=-\sum\limits_{k\neq 0} \fl{\widehat{f}_k}{k^2} e^{ikx}.
\]
Let $R=R(v,t,s)$ be a term that depends on the function values $v(t+\xi)$ for $0\le \xi\le s$. We say that $R\in \mathcal{R}_\beta(s^\alpha)$ if and only if
\[\|R(v,t,s)\|_r\le C s^\alpha,\]
where $C$ depends on $\sup\limits_{0\le\xi\le s}\|v(t+\xi)\|_{r+\beta}$. For simplicity, we write $f=g+\mathcal{R}_\beta(s^\alpha)$ whenever $f=g+R$ with $R\in \mathcal{R}_\beta(s^\alpha)$.

Next we present some lemmas which will be used frequently afterwards.
\begin{lemma}\label{eix}
For all $x, y\in\mathbb{R}$ and $0\le \alpha\le 1$, it holds that
$$|e^{ix}-1|\le 2^{1-\alpha}|x|^\alpha,\quad
|e^{ix}-1-ix|\le 2^{1-2\alpha} |x|^{1+\alpha},\quad |e^{i(x+y)}+1-e^{ix}-e^{iy}|\le 2^{2-2\alpha}|x|^\alpha |y|^\alpha.$$
\end{lemma}
\emph{Proof.} The first assertion follows by combining the following inequalities
\[|e^{ix}-1|\le 2,\qquad |e^{ix}-1|\le |x|.\]
For the second expansion, we note that
\[|e^{ix}-1-ix|\le x^2/2,\qquad |e^{ix}-1-ix|\le 2|x|,\]
which directly gives the result by
$|e^{ix}-1-ix|\le (x^2/2)^\alpha (2|x|)^{1-\alpha}=2^{1-2\alpha}|x|^{1+\alpha}$.
The last assertion follows from the first bound by noting that
\[|e^{i(x+y)}+1-e^{ix}-e^{iy}|=|e^{ix}-1||e^{iy}-1|\le 2^{2-2\alpha}|x|^\alpha |y|^\alpha,\]
which completes the proof.
\hfill $\square$ \bigskip

\begin{lemma}\label{varp12}
For all $t\in \mathbb{R}$, $\gamma\ge 0$ and $f\in H^\gamma$, we have
\[\|\psi_1(it\p_x^2)f\|_\gamma\le \|f\|_\gamma,\quad \|\psi_2(it\p_x^2)f\|_\gamma\le \|f\|_\gamma/2.\]
\end{lemma}
\emph{Proof.} For $z\in \mathbb{R}$, it can be easily obtained that
\[|\psi_1(iz)|\le 1,\quad |\psi_2(iz)|\le \int_0^1  sds\le \fl{1}{2}.\]
This directly gives
\[\|\psi_1(it\p_x^2)f\|_\gamma^2=\sum\limits_{k\in\mathbb{Z}}(1+k^2)^\gamma
|\psi_1(-itk^2)|^2|\widehat{f}_k|^2\le \|f\|_\gamma^2,\]
\[\|\psi_2(it\p_x^2)f\|_\gamma^2=\sum\limits_{k\in\mathbb{Z}}(1+k^2)^\gamma
|\psi_2(-itk^2)|^2|\widehat{f}_k|^2\le \|f\|_\gamma^2/4,\]
which completes the proof.
\hfill $\square$ \bigskip

\begin{lemma}\label{Ap}
For all $t\in \mathbb{R}$, $\gamma\ge 0$ and $f\in H^\gamma$, we have
\begin{align*}
&\|\langle \p_x^2\rangle_c^{-1} f\|_\gamma\le \|f\|_\gamma/\sqrt{c},\quad
 \|\p_x^2 \langle\p_x^2\rangle_c^{-1} f\|_\gamma\le \|f\|_\gamma,\\
&\|e^{itA}f\|_\gamma=\|f\|_\gamma,\quad  \|Af\|_\gamma\le C_1\|f\|_\gamma,\quad \|(e^{itA}-1)f\|_\gamma\le C_1 t\|f\|_\gamma,
 \end{align*}
where $C_1=\max\{1,c\}$.
\end{lemma}
\emph{Proof.}
The first two assertions are obvious by noting that the Fourier factors satisfy
\[\fl{1}{\sqrt{c+k^2+k^4}}\le \fl{1}{\sqrt{c}},\quad \left|\fl{-k^2}{\sqrt{c+k^2+k^4}}\right|\le 1.\]
In view of the fact that $A$ acts as the Fourier multiplier $A_k=\sqrt{c+k^2+k^4}-k^2$, it is bounded for $k\neq 0$ as
\be\label{Ak}
A_k=\sqrt{c+k^2+k^4}-k^2=k^2(\sqrt{1+1/k^2+c/k^4}-1)\le \fl{1}{2}+\fl{c}{2k^2}\le C_1,\ee
where $C_1=\max\{1,c\}$. Here we have used the bound $\sqrt{1+x}\le 1+x/2$ for $x\ge 0$. Hence we get
\begin{align*}
\|Af\|_\gamma^2&=\sum\limits_{k\in \mathbb{Z}}(1+k^2)^\gamma A_k^2|\widehat{f}_k|^2
\le c|\widehat{f}_0|^2+C_1^2\sum\limits_{k\neq 0}(1+k^2)^\gamma|\widehat{f}_k|^2\\
&\le C_1^2\sum\limits_{k\in \mathbb{Z}}(1+k^2)^\gamma|\widehat{f}_k|^2=C_1^2\|f\|_\gamma^2.
\end{align*}
This together with the property that $|e^{ix}-1|\le |x|$ gives $\|(e^{itA}-1)f\|_\gamma\le C_1t\|f\|_\gamma$, which completes the proof.
\hfill $\square$ \bigskip

\section{A first-order exponential-type integrator}
In this section we derive a first-order exponential-type integration scheme for the solution of equation \eqref{NF}. The construction is based on the Duhamel's formula \eqref{Duhm} and a first-order approximation.

First, by applying Lemma \ref{Ap}, \eqref{Duhm} and the bilinear estimate \eqref{bi},
we can get the first-order approximation $w(t_n+s)\approx w(t_n)$ for $|s|\le \tau$.

\begin{lemma}
For $r>1/2$, we have
\be\label{app1}
\|w(t_n+s)-w(t_n)\|_r\le 2\,C_1s\sup\limits_{0\le y\le s}
\|w(t_n+y)\|_r+D_r s\sup\limits_{0\le y\le s}
\|w(t_n+y)\|_r^2.
\ee
\end{lemma}
\emph{Proof.} Thanks to the fact that $\|\langle \p_x^2\rangle_c^{-1}f\|_r\le \fl{1}{\sqrt{c}}\|f\|_r$, it follows from \eqref{Duhm} and Lemma \ref{Ap} that
\begin{align*}
\|w(t_n+s)-w(t_n)\|_r&\le \|(e^{isA}-1)w(t_n)\|_r+\sqrt{c}\,s\sup\limits_{0\le y\le s}
\|w(t_n+y)\|_r\\
&\quad+\fl{s}{4}\sup\limits_{0\le y\le s}\left\|\left(e^{-i(t_n+y)\p_x^2}w(t_n+y)+e^{i(t_n+y)\p_x^2}\overline{w}(t_n+y)\right)^2\right\|_r\\
&\le C_1s\|w(t_n)\|_r+\sqrt{c}\,s\sup\limits_{0\le y\le s}
\|w(t_n+y)\|_r+D_r s\sup\limits_{0\le y\le s}
\|w(t_n+y)\|_r^2\\
&\le 2\,C_1s \sup\limits_{0\le y\le s}
\|w(t_n+y)\|_r+D_r s\sup\limits_{0\le y\le s}
\|w(t_n+y)\|_r^2,
\end{align*}
which completes the proof.
\hfill $\square$ \bigskip

With this approximation and Lemma \ref{Ap}, we can rewrite \eqref{Duhm} as
\begin{align*}
w(t_n+\tau)&=e^{i\tau A}w(t_n)-\fl{ic\tau}{2}\langle \p_x^2\rangle_c^{-1}
e^{i\tau A}w(t_n)-\fl{ic}{2}\langle \p_x^2\rangle_c^{-1}e^{i\tau A}\int_0^\tau e^{2i(t_n+s)\p_x^2}\overline{w(t_n)}ds\\
&\,\,-\fl{i\p_x^2}{4}\langle \p_x^2\rangle_c^{-1}e^{i\tau A}\int_0^\tau e^{i(t_n+s)\p_x^2}\left[e^{-i(t_n+s)\p_x^2}w(t_n)+e^{i(t_n+s)\p_x^2}\overline{w(t_n)}\right]^2ds
+\mathcal{R}_0(\tau^2).
\end{align*}
Twisting the variable back, we get
\begin{align}
u(t_n+\tau)&=e^{-i(t_n+\tau)\p_x^2}w(t_n+\tau)\nn\\
&=e^{i\tau \langle\p_x^2\rangle_c}u(t_n)-\fl{ic\tau}{2}\langle \p_x^2\rangle_c^{-1}
e^{i\tau \langle \p_x^2\rangle_c}u(t_n)-\fl{ic}{2}\langle \p_x^2\rangle_c^{-1}e^{i\tau \langle \p_x^2\rangle_c}
\int_0^\tau e^{2is\p_x^2}\overline{u(t_n)}ds\nn\\
&\quad-\fl{i\p_x^2}{4}\langle \p_x^2\rangle_c^{-1}e^{i\tau \langle\p_x^2\rangle_c}\int_0^\tau e^{is\p_x^2}\left[e^{-is\p_x^2}u(t_n)+e^{is\p_x^2}\overline{u(t_n)}\right]^2ds
+\mathcal{R}_0(\tau^2)\nn\\
&=e^{i\tau \langle\p_x^2\rangle_c}u(t_n)-\fl{ic\tau}{2}\langle \p_x^2\rangle_c^{-1}
e^{i\tau \langle \p_x^2\rangle_c}u(t_n)-\fl{ic\tau}{2}\langle \p_x^2\rangle_c^{-1}e^{i\tau \langle \p_x^2\rangle_c}\psi_1(2i\tau\p_x^2)\overline{u(t_n)}\nn\\
&\quad-\fl{i\p_x^2}{4}\langle \p_x^2\rangle_c^{-1}e^{i\tau \langle \p_x^2\rangle_c}\big[I_1^\tau(u(t_n))+2I_2^\tau(u(t_n))+
I_0^\tau(u(t_n))\big]+\mathcal{R}_0(\tau^2),\label{approx}
\end{align}
where
\be\label{I120}
I_1^\tau(f)=\int_0^\tau e^{is\p_x^2}(e^{-is\p_x^2}f)^2ds,\,\,
I_2^\tau(f)=\int_0^\tau e^{is\p_x^2}|e^{-is\p_x^2}f|^2ds,\,\,
I_0^\tau(f)=\int_0^\tau e^{is\p_x^2}(e^{is\p_x^2}\overline{f})^2ds.
\ee
The integral in $I_1$ can be expressed in terms of the Fourier coefficients as follows
\begin{align}
I_1^\tau(f)&=\sum\limits_{k_1,k_2\in\mathbb{Z}}
\int_0^\tau e^{is(k_1^2+k_2^2-(k_1+k_2)^2)}\widehat{f}_{k_1}\widehat{f}_{k_2}ds\, e^{i(k_1+k_2)x}\nn\\
&=\sum\limits_{k_1\neq 0, k_2\neq 0}\fl{1-e^{-2i\tau k_1k_2}}{2ik_1 k_2}\widehat{f}_{k_1}\widehat{f}_{k_2}e^{i(k_1+k_2)x}+2\tau \widehat{f}_0\sum\limits_{k\in\mathbb{Z}} \widehat{f}_ke^{ikx}-\tau\widehat{f}_0^2\nn\\
&=\fl{i}{2}\left[(\partial_x^{-1}f)^2-
e^{i\tau \p_x^2}(e^{-i\tau\p_x^2}\partial_x^{-1}f)^2\right]+2\tau \widehat{f}_0 f-\tau\widehat{f}_0^2.\label{I1d}
\end{align}
A similar calculation yields that
\begin{align}
&\hspace{-4mm}I_2^\tau(f)=\sum\limits_{k_1,k_2\in\mathbb{Z}}
\int_0^\tau e^{is(k_1^2-k_2^2-(k_1-k_2)^2)}\widehat{f}_{k_1}\overline{\widehat{f}_{k_2}}ds\, e^{i(k_1-k_2)x}\nn\\
&=\sum\limits_{k_1\neq k_2, k_2\neq 0}
\fl{e^{i\tau(k_1^2-k_2^2-(k_1-k_2)^2)}-1}{2i k_2(k_1-k_2)}\widehat{f}_{k_1}\overline{\widehat{f}_{k_2}}e^{i(k_1-k_2)x}
+\tau \overline{\widehat{f}_0}\sum\limits_{k\in\mathbb{Z}} \widehat{f}_ke^{ikx}+\tau\|f\|^2-\tau|\widehat{f}_0|^2\nn\\
&=-\fl{i}{2}\partial_x^{-1}e^{i\tau\p_x^2}\left[(e^{-i\tau\p_x^2}f)(
e^{i\tau\p_x^2}\partial_x^{-1}\overline{f})\right]+\fl{i}{2}\partial_x^{-1}\left[f(
\partial_x^{-1}\overline{f})\right]
+\tau \overline{\widehat{f}_0} f+\tau\|f\|^2-\tau|\widehat{f}_0|^2.\label{I2d}
\end{align}
For $I_0^\tau(f)$, a direct computation gives that
\begin{align*}
I_0^\tau(f)&=\sum\limits_{k_1,k_2\in\mathbb{Z}}
\overline{\widehat{f}_{k_1}}\,\overline{\widehat{f}_{k_2}} e^{-i(k_1+k_2)x}
\int_0^\tau e^{-is\left[2(k_1+k_2)^2-2k_1k_2\right]} ds\\
&=\sum\limits_{k_1,k_2\in\mathbb{Z}}
\overline{\widehat{f}_{k_1}}\,\overline{\widehat{f}_{k_2}} e^{-i(k_1+k_2)x}
\int_0^\tau e^{-2is(k_1+k_2)^2}(1+e^{2isk_1k_2}-1) ds\\
&=I_3^\tau(f)+P_1^\tau(f),
\end{align*}
where
\be\label{I3d}
I_3^\tau(f)=\tau\sum\limits_{k_1,k_2\in\mathbb{Z}}
\psi_1(-2i\tau(k_1+k_2)^2) \overline{\widehat{f}_{k_1}}\,\overline{\widehat{f}_{k_2}}e^{-i(k_1+k_2)x}=\tau \psi_1(2i\tau\p_x^2)(\overline{f})^2,
\ee
and by Lemma \ref{eix},
\begin{align*}
\|P_1^\tau(f)\|_r^2&=\sum\limits_{l\in\mathbb{Z}}(1+l^2)^r\Big|
\sum\limits_{k_1+k_2=l}\overline{\widehat{f}_{k_1}}\,\overline{\widehat{f}_{k_2}}
\int_0^\tau (e^{2isk_1k_2}-1)ds\Big|^2\\
&\le 4\tau^{2+2\gamma}\sum\limits_{l\in\mathbb{Z}}(1+l^2)^r\Big(\sum\limits_{k_1+ k_2=l}|k_1|^\gamma|k_2|^\gamma |\widehat{f}_{k_1}||\widehat{f}_{k_2}|\Big)^2\\
&\le 4\tau^{2+2\gamma}\sum\limits_{l\in\mathbb{Z}}(1+l^2)^r\Big(\sum\limits_{k_1+ k_2=l}(1+k_1^2)^{\gamma/2}(1+k_2^2)^{\gamma/2} |\widehat{f}_{k_1}||\widehat{f}_{k_2}|\Big)^2.
\end{align*}
Define
$$g(x)=\sum\limits_{k\in \mathbb{Z}}(1+k^2)^{\gamma/2}|\widehat{f}_k|e^{ikx}.$$
This implies $\|g\|_r=\|f\|_{r+\gamma}$ and
\be\label{r2}
\|P_1^\tau(f)\|_r\le 2\tau^{1+\gamma}\|g^2\|_r\le 2D_r\tau^{1+\gamma}\|f\|_{r+\gamma}^2.\ee

Combining the above approximation and noting that $I_1$, $I_2$ are exact, we can get
that
\begin{align}
&\hspace{-6mm}u(t_n+\tau)=e^{i\tau \langle\p_x^2\rangle_c}u(t_n)-\fl{ic\tau}{2}\langle \p_x^2\rangle_c^{-1}
e^{i\tau \langle \p_x^2\rangle_c}u(t_n)-\fl{ic\tau}{2}\langle \p_x^2\rangle_c^{-1}e^{i\tau \langle \p_x^2\rangle_c}\psi_1(2i\tau\p_x^2)\overline{u(t_n)}\nn\\
&-\fl{i\p_x^2}{4}\langle \p_x^2\rangle_c^{-1}e^{i\tau \langle \p_x^2\rangle_c}\big[I_1^\tau(u(t_n))+2I_2^\tau(u(t_n))+
I_3^\tau(u(t_n))\big]+\mathcal{R}_0(\tau^2)+\mathcal{R}_\gamma(\tau^{1+\gamma}).\label{approx}
\end{align}
Let $u^0=u_0$. The above calculation motivates us to define the scheme for $n\ge 0$ as
\be\label{unp}
u^{n+1}=\Phi^\tau (u^n),
\ee
where
\be\label{schu}
\begin{split}
\Phi^\tau(f)&=e^{i\tau\langle\p_x^2\rangle_c }f-\fl{ic\tau}{2}\langle \p_x^2\rangle_c^{-1}
e^{i\tau\langle\p_x^2\rangle_c }f-\fl{ic\tau}{2}\langle \p_x^2\rangle_c^{-1}e^{i\tau \langle\p_x^2\rangle_c}\psi_1(2i\tau\p_x^2)\overline{f}\\
&\quad-\fl{i\p_x^2}{4}\langle \p_x^2\rangle_c^{-1}e^{i\tau \langle \p_x^2\rangle_c}\big[I_1^\tau(f)+2I_2^\tau(f)+I_3^\tau(f)\big].
\end{split}
\ee
The terms $I_1^\tau$, $I_2^\tau$ and $I_3^\tau$ are defined in \eqref{I1d}, \eqref{I2d} and \eqref{I3d}, respectively.

Regarding the stability analysis, we have the following estimate.
\begin{lemma}\label{stab1}
Let $r>1/2$ and $f, g\in H^{r}$. Then for all $\tau\ge 0$ we have
\[\|\Phi^\tau(f)-\Phi^\tau(g)\|_r\le (1+M\tau)\|f-g\|_r,\]
where $M=\sqrt{c}+D_r(\|f\|_r+\|g\|_r)$.
\end{lemma}
\emph{Proof.} For simplicity of notation, let $Q=\|f\|_r+\|g\|_r$. Using \eqref{I120} and \eqref{bi}, we have
\be\label{I1e}
\begin{split}
\|I_1^\tau(f)-I_1^\tau(g)\|_r&\le\tau\sup\limits_{0\le s\le\tau}
\|(e^{-is\p_x^2}f)^2-(e^{-is\p_x^2}g)^2\|_r\\
&\le D_{r}\tau\sup\limits_{0\le s\le\tau}
\|e^{-is\p_x^2}(f+g)\|_r\|e^{-is\p_x^2}(f-g)\|_r\le Q D_r\tau\|f-g\|_r.
\end{split}
\ee
A similar calculation gives that
\be\label{I2e}
\|I_2^\tau(f)-I_2^\tau(g)\|_r\le Q D_r\tau\|f-g\|_r.
\ee
With the help of Lemma \ref{varp12}, we obtain
\be\label{I3e}
\|I_3^\tau(f)-I_3^\tau(g)\|_r\le \tau \|f^2-g^2\|_r\le Q D_r\tau\|f-g\|_r.
\ee
Combining \eqref{I1e}-\eqref{I3e} and applying Lemmas \ref{varp12} and \ref{Ap}, we finally derive
\begin{align*}
\|\Phi^\tau(f)-\Phi^\tau(g)\|_r&\le (1+\sqrt{c}\tau)\|f-g\|_r+Q D_r\tau\|f-g\|_r\\
&\le \left[1+\tau(\sqrt{c}+Q D_r)\right]\|f-g\|_r,
\end{align*}
which completes the proof.
\hfill $\square$ \bigskip

The local error estimate \eqref{approx} and the stability property allow us to prove the following error bound.

\begin{theorem}\label{first-th}
Let $r>1/2$ and $0<\gamma\le 1$. Assume that the exact solution of \eqref{NF} satisfies $u\in H^{r+\gamma}$ for $0\le t\le T$. Then there exists a constant $\tau_0>0$ such that for all step sizes
$0<\tau\le \tau_0$ and $t_n=n\tau\le T$ we have that the global error of \eqref{znp} is bounded by
\[\|u(t_n)-u^n\|_r\le C\tau^\gamma,\]
where $C$ depends on $T$, $c$, $r$ and $\|u\|_{L^\infty(0,T;H^{r+\gamma})}$.
\end{theorem}
\emph{Proof.} It follows from \eqref{approx} that
\[\|u(t_{k+1})-\Phi^\tau(u(t_k))\|_r\le M_1\tau^{1+\gamma},\]
where $M_1$ depends on $\|u\|_{L^\infty(0,T; H^{r+\gamma})}$.
The triangle inequality yields
\begin{align*}
\|u(t_{k+1})-u^{k+1}\|_r&=\|u(t_{k+1})-\Phi^\tau(u^k)\|_r\\
&\le \|u(t_{k+1})-\Phi^\tau(u(t_k))\|_r+\|\Phi^\tau(u(t_k))-\Phi^\tau(u^k)\|_r\\
&\le M_1\tau^{1+\gamma}+\|\Phi^\tau(u(t_k))-\Phi^\tau(u^k)\|_r.
\end{align*}
By applying Lemma \ref{stab1} for $u^k\in H^r (0\le k\le n)$, we get
\begin{align*}
&\hspace{-6mm}\|u(t_{n+1})-u^{n+1}\|_r\le M_1\tau^{1+\gamma}+e^{\tau L}\|u(t_n)-u^n\|_r\\
&\le M_1\tau^{1+\gamma}+e^{\tau L}\left(M_1\tau^{1+\gamma}+e^{\tau L}\|u(t_{n-1})-u^{n-1}\|_r\right)\\
&\le M_1\tau^{1+\gamma}\sum\limits_{k=0}^ne^{t_k L}\le \fl{M_1}{L}e^{TL}\tau^\gamma ,
\end{align*}
where $L$ depends on $\sup\limits_{0\le k\le n}\|u(t_k)\|_r$ and $\sup\limits_{0\le k\le n}\|u^k\|_r$. Then the assertion follows by induction, respectively, a \emph{Lady Windermere's fan} argument (cf. for example \cite{holden2013operator, lubich2008splitting}).
\hfill $\square$ \bigskip

For the scheme $u^{n+1}$ defined in \eqref{unp}, we set
\be\label{znp}
z^{n+1}=\fl{1}{2}(u^{n+1}+\overline{u^{n+1}}),\quad z_t^{n+1}=\fl{i}{2}\langle \p_x^2\rangle_c (u^{n+1}-\overline{u^{n+1}}).
\ee
In view of the relationship between $u$ and $z$, Theorem \ref{first-th} allows us to get the convergence of $z^n$.
\begin{corollary}\label{1-th}
Let $r>1/2$ and $0<\gamma\le 1$. Assume that the exact solution of \eqref{GBp} satisfies
\[z\in H^{r+\gamma},\quad z_t\in H^{r+\gamma-2},\]
which implies that the solution of \eqref{NF} satisfies $u\in H^{r+\gamma}$
for $0\le t\le T$.
Then there exists a constant $\tau_0>0$ such that for all step sizes
$0<\tau\le \tau_0$ and $t_n=n\tau\le T$ we have that the global error of \eqref{znp} is bounded by
\[\|z(t_n)-z^n\|_r+\|z_t(t_n)-z_t^n\|_{r-2}\le C\tau^\gamma,\]
where $C$ depends on $T$, $c$, $r$ and $\|z\|_{L^\infty(0,T;H^{r+\gamma})}+\|z_t\|_{L^\infty(0,T;H^{r+\gamma-2})}$.
\end{corollary}

\begin{remark}
Due to the well-posedness result from \cite{kishimoto2013sharp}, the GB equation is locally well-posed for initial data $z(0)\in H^s$, $z_t(0)\in H^{s-2}$ with $s\ge -1/2$. Thus the regularity required in Corollary \ref{1-th} can be obtained whenever the initial data has the same smoothness.
\end{remark}
\section{A second-order exponential-type integrator}
To construct the second-order scheme, we need the following expansions of the operators
$e^{is A}$ and $e^{is \p_x^2}$ based on the bounds in Lemma \ref{eix}.
\begin{lemma}\label{expan}
For all $s\in\mathbb{R}$ and $v\in H^r(r\ge 0)$, we have
\be\label{A2c}
e^{\pm is A} v=v\pm is Av+\mathcal{R}_0(s^2).
\ee
If $0\le \gamma \le 1$ and $v\in H^{r+2\gamma}$, then
\be\label{L1c}
e^{\pm is \partial_x^2}v=v+\mathcal{R}_{2\gamma}(s^\gamma).
\ee
\end{lemma}
\emph{Proof.} The first assertion \eqref{A2c} follows from the second-order expansion in Lemma \ref{eix} and the boundedness of the operator $A$ (cf. Lemma \ref{Ap}). Similarly, the first-order expansion in Lemma \ref{eix} yields the second identity.
\hfill $\square$ \bigskip

Using \eqref{Duhm} and \eqref{app1}, we can approximate $w(t_n+s)$ as follows:
\begin{align*}
w(t_n+s)&=e^{is A}w(t_n)-\fl{ic}{2}\langle \p_x^2\rangle_c^{-1}\int_0^s
e^{i(s-y)A}\left[w(t_n)+e^{2i(t_n+y)\p_x^2}\overline{w(t_n)}\right]dy\\
&\quad-\fl{i\p_x^2}{4}\langle \p_x^2\rangle_c^{-1}\int_0^s e^{i(s-y)A}e^{i(t_n+y)\p_x^2}\left[e^{-iy\p_x^2}u(t_n)+e^{iy\p_x^2}
\overline{u(t_n)}\right]^2dy+\mathcal{R}_0(s^2),
\end{align*}
where we used the identity $w(t_n)=e^{it_n\p_x^2}u(t_n)$. Let $0<\gamma\le 1$. Using \eqref{A2c} and \eqref{L1c}, we get
\begin{align}
w(t_n+s)&=e^{is A}w(t_n)-\fl{ics}{2}\langle \p_x^2\rangle_c^{-1}\left[w(t_n)+e^{2it_n\p_x^2}\overline{w(t_n)}\right]\nn\\
&\quad-\fl{is\p_x^2}{4}\langle \p_x^2\rangle_c^{-1} e^{it_n\p_x^2}\left[u(t_n)+
\overline{u(t_n)}\right]^2+\mathcal{R}_0(s^2)+\mathcal{R}_{2\gamma}(s^{1+\gamma})\nn\\
&=w(t_n)+isv(t_n)+\mathcal{R}_{2\gamma}(s^{1+\gamma}),\label{wap}
\end{align}
where
\be\label{vd}
v(t_n)=Aw(t_n)-\fl{1}{4}\langle \p_x^2\rangle_c^{-1}\left[2c\big(w(t_n)+ e^{2it_n\p_x^2}\overline{w(t_n)}\big)+\p_x^2 e^{it_n\p_x^2}(u(t_n)+
\overline{u(t_n)})^2\right].
\ee
Plugging \eqref{wap} into \eqref{Duhm} and applying Lemma \ref{expan}, we get
\begin{align}
&\hspace{-2mm}w(t_n+\tau)=e^{i\tau A}w(t_n)+\mathcal{R}_{2\gamma}(\tau^{2+\gamma})\nn\\
&\quad\,\,-\fl{ic}{2}\langle \p_x^2\rangle_c^{-1}e^{i\tau A}\int_0^\tau
e^{-isA}\big[w(t_n)+isv(t_n)+e^{2i(t_n+s)\p_x^2}(\overline{w(t_n)}-is\overline{v(t_n)})
\big]ds\nn\\
&\quad\,\,-\fl{i\p_x^2}{4}\langle \p_x^2\rangle_c^{-1}e^{i\tau A} e^{it_n \p_x^2}\int_0^\tau e^{is(\p_x^2-A)}\big[e^{-is\p_x^2}(u(t_n)+is\mu(t_n))+e^{is\p_x^2}
(\overline{u(t_n)}-is\overline{\mu(t_n)})\big]^2ds\nn\\
&\,\,\,=e^{i\tau A}w(t_n)-\fl{ic}{2}\langle \p_x^2\rangle_c^{-1}e^{i\tau A}\Big[\int_0^\tau \left[(1-isA)w(t_n)+isv(t_n)\right]ds\Big.\nn\\ \Big.
&\qquad\qquad\qquad\qquad\qquad\qquad\qquad\Big.+e^{2it_n\p_x^2}\int_0^\tau e^{2is\p_x^2}\big[(1-isA)\overline{w(t_n)}-is\overline{v(t_n)}\big]ds\Big]\nn\\
&\quad\,\,\,\,-\fl{i\p_x^2}{4}\langle\p_x^2\rangle_c^{-1}e^{i\tau A} e^{it_n \p_x^2}\int_0^\tau
e^{is\p_x^2}\left[(1-isA)\big(e^{-is\p_x^2}u(t_n)+e^{is\p_x^2}\overline{u(t_n)}\big)^2
\right.\nn\\
&\quad\,\,\,\,\left.+2is \big(e^{-is\p_x^2}u(t_n)+e^{is\p_x^2}\overline{u(t_n)}\big)
\big(e^{-is\p_x^2}\mu(t_n)-e^{is\p_x^2}\overline{\mu(t_n)}\big)\right]ds+
\mathcal{R}_{2\gamma}(\tau^{2+\gamma}),
\end{align}
where
\[\mu(t_n)=e^{-it_n\p_x^2}v(t_n)=\mathcal{P}(u(t_n)),\]
with
\be\label{Pd}
\mathcal{P}(f)=Af-\fl{1}{4}\langle \p_x^2\rangle_c^{-1}\left[2c\big(f+\overline{f}\big)+\p_x^2 (f+\overline{f})^2\right].
\ee
Twisting the variable back, we get the approximation for $u(t_n+\tau)$ as
\begin{align}
u(t_n+\tau)&=e^{i\tau \langle \p_x^2\rangle_c}u(t_n)-\fl{ic}{2}\langle \p_x^2\rangle_c^{-1}e^{i\tau \langle\p_x^2\rangle_c}\Big[\int_0^\tau \left[(1-isA)u(t_n)+is\mu(t_n)\right]ds\Big.\nn\\ \Big.
&\qquad\qquad\qquad\qquad\qquad\qquad\qquad\quad\Big.+\int_0^\tau e^{2is\p_x^2}\big[(1-isA)\overline{u(t_n)}-is\overline{\mu(t_n)}\big]ds\Big]\nn\\
&-\fl{i\p_x^2}{4}\langle\p_x^2\rangle_c^{-1}e^{i\tau \langle\p_x^2\rangle_c} \int_0^\tau
e^{is\p_x^2}\left[(1-isA)\big(e^{-is\p_x^2}u(t_n)+e^{is\p_x^2}\overline{u(t_n)}\big)^2
\right.\nn\\
&\left.+2is \big(e^{-is\p_x^2}u(t_n)+e^{is\p_x^2}\overline{u(t_n)}\big)
\big(e^{-is\p_x^2}\mu(t_n)-e^{is\p_x^2}\overline{\mu(t_n)}\big)\right]ds+\mathcal{R}_0(\tau^3)+
\mathcal{R}_{2\gamma}(\tau^{2+\gamma})\nn\\
&\hspace{-5mm}=e^{i\tau \langle\p_x^2\rangle_c} u(t_n)+\mathcal{L}(u(t_n))+\fl{\p_x^2}{4}\langle\p_x^2\rangle_c^{-1}e^{i\tau \langle\p_x^2\rangle_c}\Big(-i\big[I_1^\tau(u(t_n))+2I_2^\tau(u(t_n))+I_0^\tau(u(t_n))\big]\Big.\nn\\
&-A\big[J_1^\tau(u(t_n),u(t_n))+2J_2^\tau(u(t_n),u(t_n))+J_0^\tau(u(t_n),u(t_n))\big]
+2J_1^\tau(u(t_n),\mu(t_n))\nn\\
&\Big.+2\left[J_2^\tau(\mu(t_n),u(t_n))-
J_0^\tau(u(t_n),\mu(t_n))-J_2^\tau(u(t_n),\mu(t_n))\right]\Big)
+\mathcal{R}_{2\gamma}(\tau^{2+\gamma}),\label{unapp1}
\end{align}
where $I_1^\tau$, $I_2^\tau$, $I_0^\tau$ are defined as \eqref{I120},
\be\label{Ld}
\mathcal{L}(f)=-\fl{ic\tau}{2}\langle \p_x^2\rangle_c^{-1}e^{i\tau \langle\p_x^2 \rangle_c}\Big[f-\fl{i\tau}{2}\left(Af-\mathcal{P}(f)\right)+\psi_1(2i\tau\p_x^2)
\overline{f}-i\tau\psi_2(2i\tau\p_x^2)\big(A\overline{f}+\overline{\mathcal{P}(f)}\big)\Big],
\ee
and
\be\label{Jd}
\begin{split}
&J_1^\tau(f,g)=\int_0^\tau
se^{is\p_x^2}\big(e^{-is\p_x^2}f\big)\big(e^{-is\p_x^2}g\big)ds,\quad
J_2^\tau(f,g)=\int_0^\tau
se^{is\p_x^2}\big(e^{-is\p_x^2}f\big)\big(e^{is\p_x^2}\overline{g}\big) ds,\\
&J_0^\tau(f,g)=\int_0^\tau
se^{is\p_x^2}\big(e^{is\p_x^2}\overline{f}\big)\big(e^{is\p_x^2}\overline{g}\big)ds.
\end{split}
\ee

The integral in $J_1^\tau(f,g)$ can be expressed in terms of the Fourier coefficients as follows
\begin{align}
J_1^\tau(f,g)&=\sum\limits_{k_1,k_2\in\mathbb{Z}}
\int_0^\tau s e^{is(k_1^2+k_2^2-(k_1+k_2)^2)}\widehat{f}_{k_1}\widehat{g}_{k_2}ds\, e^{i(k_1+k_2)x}\nn\\
&=\sum\limits_{k_1\neq 0, k_2\neq 0}\left[\fl{i\tau e^{-2i\tau k_1k_2}}{2k_1 k_2}+\fl{e^{-2i\tau k_1k_2}-1}{4k_1^2k_2^2}\right]\widehat{f}_{k_1}\widehat{g}_{k_2}e^{i(k_1+k_2)x}\nn\\
&\quad+\fl{\tau^2}{2} \widehat{f}_0\sum\limits_{k\in\mathbb{Z}} \widehat{g}_ke^{ikx}+\fl{\tau^2}{2} \widehat{g} _0\sum\limits_{k\in\mathbb{Z}} \widehat{f}_ke^{ikx}-\fl{\tau^2}{2}\widehat{f}_0\widehat{g}_0\nn\\
&=-\fl{i\tau}{2}e^{i\tau \p_x^2}\left[\big(\partial_x^{-1}e^{-i\tau \p_x^2}f\big)
\big(\partial_x^{-1}e^{-i\tau \p_x^2}g\big)\right]-\fl{1}{4}\big(\partial_x^{-2}f\big)
\big(\partial_x^{-2}g\big)\nn\\
&\quad+\fl{1}{4}e^{i\tau \p_x^2}\left[\big(\partial_x^{-2}e^{-i\tau \p_x^2}f\big)
\big(\partial_x^{-2}e^{-i\tau \p_x^2}g\big)\right]
+\fl{\tau^2}{2} \widehat{f}_0g+\fl{\tau^2}{2} \widehat{g} _0f-\fl{\tau^2}{2}\widehat{f}_0\widehat{g}_0.\label{J1d}
\end{align}
A similar calculation yields that
\begin{align}
&\hspace{-4mm}J_2^\tau(f,g)=\sum\limits_{k_1,k_2\in\mathbb{Z}}
\int_0^\tau s e^{is(k_1^2-k_2^2-(k_1-k_2)^2)}\widehat{f}_{k_1}\overline{\widehat{g}_{k_2}}ds\, e^{i(k_1-k_2)x}\nn\\
&=\sum\limits_{k_1\neq k_2, k_2\neq 0}\left[\fl{\tau e^{2i\tau k_2(k_1-k_2)}}{2ik_2 (k_1-k_2)}+\fl{e^{2i\tau k_2(k_1-k_2)}-1}{4k_2^2(k_1-k_2)^2}\right]\widehat{f}_{k_1}\overline{\widehat{g}_{k_2}}
e^{i(k_1-k_2)x}\nn\\
&\quad+\fl{\tau^2}{2}\overline{\widehat{g}_0}\sum\limits_{k} \widehat{f}_ke^{ikx}+\fl{\tau^2}{2}\sum\limits_{k} \widehat{f}_k\overline{\widehat{g}_{k}}-\fl{\tau^2}{2}\widehat{f}_0
\overline{\widehat{g}_{0}}\nn\\
&=\fl{i\tau}{2}\partial_x^{-1}e^{i\tau\p_x^2}\left[\big(e^{-i\tau\p_x^2}f\big)\big(
e^{i\tau\p_x^2}\partial_x^{-1}\overline{g}\big)\right]+\fl{1}{4}\partial_x^{-2}
e^{i\tau\p_x^2}\left[\big(e^{-i\tau\p_x^2}f\big)\big(
e^{i\tau\p_x^2}\partial_x^{-2}\overline{g}\big)\right]\nn\\
&\quad-\fl{1}{4}\partial_x^{-2}
\left[f\big(\partial_x^{-2}\overline{g}\big)\right]
+\fl{\tau^2}{2} \overline{\widehat{g}_0} f+\fl{\tau^2}{2}(f, g)-\fl{\tau^2}{2}\widehat{f}_0
\overline{\widehat{g}_{0}},\label{J2d}
\end{align}
where $(\cdot, \cdot)$ represents the inner product in $L^2$ defined as $(f, g)=\fl{1}{2\pi}\int_\Omega f(x)\overline{g(x)}dx$.
Similarly, one obtains
\begin{align}
J_0^\tau(f,g)&=\sum\limits_{k_1,k_2\in\mathbb{Z}}
\int_0^\tau s e^{-is\left[2(k_1+k_2)^2-2k_1k_2\right]}ds\, \overline{\widehat{f}_{k_1}}\overline{\widehat{g}_{k_2}}e^{-i(k_1+k_2)x}\nn\\
&=\sum\limits_{k_1,k_2\in\mathbb{Z}}
\int_0^\tau s e^{-2is(k_1+k_2)^2}ds\, \overline{\widehat{f}_{k_1}}\overline{\widehat{g}_{k_2}}e^{-i(k_1+k_2)x}+P_0^\tau(f,g)\nn\\
&=\sum\limits_{k_1+k_2\neq 0}\left[\fl{i\tau e^{-2i\tau (k_1+k_2)^2}}{2(k_1+k_2)^2}+\fl{e^{-2i\tau (k_1+k_2)^2}-1}{4(k_1+k_2)^4}\right]\overline{\widehat{f}_{k_1}}\overline{\widehat{g}_{k_2}}e^{-i(k_1+k_2)x}\nn\\
&\quad+\fl{\tau^2}{2}\sum\limits_{k\in\mathbb{Z}}\overline{\widehat{f}_{k}}
\overline{\widehat{g}_{-k}}+P_0^\tau(f,g)=J_{3}^\tau(f,g)+P_0^\tau(f,g),\label{J0}
\end{align}
where
\be\label{J3d}
J_{3}^\tau(f,g)=-\fl{i\tau}{2}e^{2i\tau\p_x^2}\p_x^{-2}\left(\overline{fg}\right)+
\fl{i\tau}{2}\psi_1(2i\tau\p_x^2)\p_x^{-2}\left(\overline{fg}\right)
+\fl{\tau^2}{2}(\overline{f},g).
\ee
By Lemma \ref{eix}, we have
\begin{align*}
\|P_0^\tau(f,g)\|_r^2&=\sum\limits_{l\in\mathbb{Z}}(1+l^2)^r\,\Big|
\sum\limits_{k_1+k_2=l}\overline{\widehat{f}_{k_1}}\overline{\widehat{g}_{k_2}}
\int_0^\tau s e^{-2isl^2}(e^{2isk_1k_2}-1)ds\Big|^2\\
&\le \tau^{4+2\gamma}\sum\limits_{l\in\mathbb{Z}}(1+l^2)^r\Big(\sum\limits_{k_1+ k_2=l}|k_1|^\gamma|k_2|^\gamma |\widehat{f}_{k_1}||\widehat{g}_{k_2}|\Big)^2\\
&\le \tau^{4+2\gamma}\sum\limits_{l\in\mathbb{Z}}(1+l^2)^r\Big(\sum\limits_{k_1+ k_2=l}(1+k_1^2)^{\gamma/2}(1+k_2^2)^{\gamma/2} |\widehat{f}_{k_1}||\widehat{g}_{k_2}|\Big)^2,
\end{align*}
which implies that
\be\label{R0}
\|P_0^\tau(f,g)\|_r\le D_r\tau^{2+\gamma}\|f\|_{r+\gamma}\|g\|_{r+\gamma}.
\ee
It follows from Lemma \ref{Ap} that
\[\|\mu(t_n)\|_{r+\gamma}\le 2C_1\|u(t_n)\|_{r+\gamma}+D_{r+\gamma}\|u(t_n)\|_{r+\gamma}^2,\]
which together with \eqref{J0} yields that
\be\label{J0app}
\begin{split}
&J_0^\tau(u(t_n),\mu(t_n))=J_{3}^\tau(u(t_n),\mu(t_n))+\mathcal{R}_\gamma(\tau^{2+\gamma}),\\
&J_0^\tau(u(t_n),u(t_n))=J_{3}^\tau(u(t_n),u(t_n))+\mathcal{R}_\gamma(\tau^{2+\gamma}).
\end{split}
\ee

Instead of approximating $I_0^\tau$ by $I_3^\tau$ with the remainder $\mathcal{R}_\gamma(\tau^{1+\gamma})$, it remains to find a refined approximation to the integral $I_0^\tau$. For this aim, we will employ the following decomposition
\[e^{-is\left[2(k_1+k_2)^2-2k_1k_2\right]}=e^{-2is(k_1+k_2)^2}
e^{2isk_1k_2}=e^{-2is(k_1+k_2)^2}+e^{2isk_1k_2}-1+
P_2^s(k_1,k_2),\]
where, by Lemma \ref{eix},
\[|P_2^s(k_1,k_2)|=|e^{-2is(k_1+k_2)^2}-1||e^{2isk_1k_2}-1|\le 4s^{1+\gamma}|k_1+k_2|^{2\gamma}|k_1||k_2|.\]
Hence we have
\begin{align*}
I_0^\tau(f)&=\sum\limits_{k_1,k_2\in\mathbb{Z}}
\overline{\widehat{f}_{k_1}}\,\overline{\widehat{f}_{k_2}} e^{-i(k_1+k_2)x}\int_0^\tau \left[e^{-2is(k_1+k_2)^2}+e^{2isk_1k_2}-1\right]ds+P_2^\tau(f),
\end{align*}
where
\begin{align*}
\|P_2^\tau(f)\|_r^2&=\sum\limits_{l\in\mathbb{Z}}(1+l^2)^r\,\Big|
\sum\limits_{k_1+k_2=l}\overline{\widehat{f}_{k_1}}\,\overline{\widehat{f}_{k_2}}
\int_0^\tau P_2^s(k_1,k_2)ds\Big|^2\\
&\le 16\tau^{4+2\gamma}\sum\limits_{l\in\mathbb{Z}}(1+l^2)^{r+2\gamma}\Big(\sum\limits_{k_1+ k_2=l}|k_1||k_2| |\widehat{f}_{k_1}||\widehat{f}_{k_2}|\Big)^2,
\end{align*}
which implies that
\[\|P_2^\tau(f)\|_r\le 4D_{r+2\gamma}\tau^{2+\gamma}\|f\|_{r+1+2\gamma}^2.\]
Hence
\begin{align}
I_0^\tau(f)&=\sum\limits_{k_1,k_2\in\mathbb{Z}}
\overline{\widehat{f}_{k_1}}\,\overline{\widehat{f}_{k_2}} e^{-i(k_1+k_2)x}\int_0^\tau \left[e^{-2is(k_1+k_2)^2}+e^{2isk_1k_2}-1\right]ds+\mathcal{R}_{1+2\gamma}(\tau^{2+\gamma})\nn\\
&=I_3^\tau(f)+\overline{I_1(f)}-\tau \overline{f}^2 +\mathcal{R}_{1+2\gamma}(\tau^{2+\gamma}),\label{I0app}
\end{align}
where $I_1^\tau$ and $I_3^\tau$ are defined as \eqref{I1d} and \eqref{I3d}, respectively.

Note that $I_1$, $I_2$, $J_1$ and $J_2$ are exact, combining \eqref{unapp1} and the approximations \eqref{J0app} and \eqref{I0app}, we get
\begin{align}
u(t_{n+1})&=e^{i\tau\langle\p_x^2\rangle_c }u(t_n)+\mathcal{L}(u(t_n))+\fl{\p_x^2}{4}\langle\p_x^2\rangle_c^{-1}e^{i\tau \langle\p_x^2\rangle_c} \Big(-i\big[I_1^\tau(u(t_n))+\overline{ I_1^\tau(u(t_n))}+2I_2^\tau(u(t_n))\big.\Big.\nn\\
&\quad\big.+I_3^\tau(u(t_n))-\tau\overline{u(t_n)}^2\big]-A\big[J_1^\tau(u(t_n),u(t_n))+2J_2^\tau(u(t_n),u(t_n))+J_3^\tau(u(t_n),u(t_n))\big]\nn\\
&\quad\Big.+2\left[J_1^\tau(u(t_n),\mu(t_n))+J_2^\tau(\mu(t_n),u(t_n))-
J_3^\tau(u(t_n),\mu(t_n))-J_2^\tau(u(t_n),\mu(t_n))\right]\Big)\nn\\
&\quad+\mathcal{R}_{1+2\gamma}(\tau^{2+\gamma}).\label{local2}
\end{align}
This motivates us to define the numerical scheme as
\be\label{unp2}
u^{n+1}=\Psi^\tau(u^n),
\ee
where
\begin{align}
\Psi^\tau(f)&=e^{i\tau\langle\p_x^2\rangle_c }f+\mathcal{L}(f)+\fl{\p_x^2}{4}\langle\p_x^2\rangle_c^{-1}e^{i\tau \langle\p_x^2\rangle_c} \Big(-i\big[I_1^\tau(f)+\overline{I_1^\tau(f)}+2I_2^\tau(f)+I_{3}^\tau(f)-\tau\overline{f}^2\big]\Big.\nn\\
&\quad-A\left[J_1^\tau(f, f)+2J_{2}^\tau(f, f)
+J_3^\tau(f, f)\right]+2\left[J_1^\tau(f,\mathcal{P}(f))+J_2^\tau(\mathcal{P}(f), f)\right]\nn\\
&\quad-2\left[J_{3}^\tau(f,\mathcal{P}(f))+J_2^\tau(f,\mathcal{P}(f))\right]\Big),\label{sch2}
\end{align}
with $\mathcal{L}$ and $\mathcal{P}$ defined in \eqref{Ld} and \eqref{Pd}, respectively,
and $I_1^\tau$, $I_2^\tau$, $I_{3}^\tau$, $J_1^\tau$, $J_{2}^\tau$, $J_{3}^\tau$ defined
in \eqref{I1d}, \eqref{I2d}, \eqref{I3d}, \eqref{J1d}, \eqref{J2d} and \eqref{J3d}, respectively.

\bigskip

Regarding the stability analysis, we have the following estimate.
\begin{lemma}
Let $r>1/2$ and $f, g\in H^{r}$. Then, for all $\tau\le 1$, we have
\be\label{stab2}
\|\Psi^\tau(f)-\Psi^\tau(g)\|_r\le (1+M\tau)\|f-g\|_r,
\ee
where $M$ depends on $c$, $r$ and $\|f\|_r+\|g\|_r$.
\end{lemma}
\emph{Proof.} We still use the notation $Q=\|f\|_r+\|g\|_r$.
By the definitions of $\mathcal{L}$ and $\mathcal{P}$, applying Lemma \ref{Ap} and \eqref{bi}, we get
\[\|\mathcal{P}(f)\|_r\le (C_1+\sqrt{c})\|f\|_r+D_r\|f\|_r^2\le (2C_1+D_r\|f\|_r)\|f\|_r,\]
\[
\|\mathcal{P}(f)-\mathcal{P}(g)\|_r\le C_1\|f-g\|_r+\sqrt{c}\|f-g\|_r+
\fl{1}{2}\left(\|f^2-g^2\|_r+\||f|^2-|g|^2\|_r\right)\le M_1\|f-g\|_r,
\]
where $M_1=2C_1+D_r Q$. This together with Lemma \ref{varp12} yields that
\begin{align*}
\|\mathcal{L}(f)-\mathcal{L}(g)\|_r&\le \fl{\tau\sqrt{c}}{2}\left[(2+C_1\tau)\|f-g\|_r+\tau\|\mathcal{P}(f)-\mathcal{P}(g)\|_r
\right]\\
&\le \tau M_2\|f-g\|_r,
\end{align*}
where $$M_2=\fl{\sqrt{c}}{2}(2+C_1\tau+M_1\tau)=\fl{\sqrt{c}}{2}(2+3C_1\tau+D_rQ\tau).$$
By the definition of $J_1$ and $J_2$, we have
\begin{align*}
\|J_1^\tau(f_1,g_1)-J_1^\tau(f_2,g_2)\|_r&\le\Big\|\int_0^\tau se^{is\p_x^2}(e^{-is\p_x^2}f_1)(e^{-is\p_x^2}(g_1-g_2))ds\Big\|_r\\
&\quad+
\Big\|\int_0^\tau se^{is\p_x^2}(e^{-is\p_x^2}g_2)(e^{-is\p_x^2}(f_1-f_2))ds\Big\|_r\\
&\le D_r\tau^2\left(\|f_1\|_r\|g_1-g_2\|_r+\|g_2\|_r\|f_1-f_2\|_r\right),
\end{align*}
which immediately gives that
\begin{align*}
\|J_1^\tau(f,f)-J_1^\tau(g,g)\|_r
&\le QD_r \tau^2\|f-g\|_r,\\
\|J_1^\tau(f,\mathcal{P}(f))-
J_1^\tau(g,\mathcal{P}(g))\|_r&\le D_r\tau^2(\|f\|_r\|\mathcal{P}(f)-\mathcal{P}(g)\|_r+\|\mathcal{P}(g)\|_r\|f-g\|_r)\\
&\le M_3\tau^2\|f-g\|_r,
\end{align*}
where
\[M_3=QD_r(2C_1+QD_r).\]
Similar calculations show that
\begin{align*}
&\|J_2^\tau(f,f)-J_2^\tau(g,g)\|_r
\le QD_r\tau^2\|f-g\|_r,\\
&\|J_2^\tau(f,\mathcal{P}(f))-
J_2^\tau(g,\mathcal{P}(g))\|_r\le M_3\tau^2\|f-g\|_r,\\
&\|J_2^\tau(\mathcal{P}(f),f)-
J_2^\tau(\mathcal{P}(g),g)\|_r\le M_3\tau^2\|f-g\|_r.
\end{align*}
For the approximation term $J_3$, we estimate the Lipschitz continuity together with the operator in front of it, which allows us to ignore the constant in $J_3$ \eqref{J3d}:
\begin{align}
\|\fl{\p_x^2}{4}\langle\p_x^2\rangle_c^{-1}e^{i\tau \langle\p_x^2\rangle_c} \left[J_{3}^\tau(f_1,g_1)-J_{3}^\tau(f_2,g_2)\right]\|_r&=
\|\fl{\p_x^2}{4}\langle\p_x^2\rangle_c^{-1} \left[J_4^\tau(f_1,g_1)-J_4^\tau(f_2,g_2)\right]\|_r\nn\\
&\le\fl{1}{4}\|J_4^\tau(f_1,g_1)-J_4^\tau(f_2,g_2)\|_r,\label{J3s}
\end{align}
where $J_4$ is obtained from $J_3$ by removing the constant
\be\label{J30d}
J_4^\tau(f,g)=-\fl{i\tau}{2}e^{2i\tau\p_x^2}\p_x^{-2}\left(\overline{fg}\right)+
\fl{i\tau}{2}\psi_1(2i\tau\p_x^2)\p_x^{-2}\left(\overline{fg}\right).
\ee
Applying Lemma \ref{varp12}, one can easily find that
\[
\|J_4^\tau(f_1,g_1)-J_4^\tau(f_2,g_2)\|_r\le D_r\tau \left(\|f_1\|_r\|g_1-g_2\|_r+\|g_2\|_r\|f_1-f_2\|_r\right),\]
which implies that
\begin{align*}
&\|J_4^\tau(f,f)-J_4^\tau(g,g)\|_r
\le QD_r\tau\|f-g\|_r,\\
&\|J_4^\tau(f,\mathcal{P}(f))-
J_4^\tau(g,\mathcal{P}(g))\|_r\le D_r\tau(\|f\|_r\|\mathcal{P}(f)-\mathcal{P}(g)\|_r+
\|\mathcal{P}(g)\|_r\|f-g\|_r)\\
&\qquad\qquad\qquad\qquad\qquad\qquad\,\,\,\,\le M_3\tau\|f-g\|_r.
\end{align*}
Combining the estimates above and noticing \eqref{J3s}, we derive that
\begin{align*}
&\hspace{-3mm}\|\Psi^\tau(f)-\Psi^\tau(g)\|_r\\
&\le \|f-g\|_r+\|\mathcal{L}(f)-\mathcal{L}(g)\|_r+ \fl{1}{4}\Big[2\sum\limits_{j=1}^2\|I_j^\tau(f)-I_j^\tau(g)\|_r+\|I_3^\tau(f)-I_3^\tau(g)\|_r\Big.\\
&\quad+\tau\|f^2-g^2\|_r+2\sum\limits_{j=1}^2\|J_j^\tau(f,\mathcal{P}(f))-
J_j^\tau(g,\mathcal{P}(g))\|_r\\
&\quad+2\left(\|J_2^\tau(\mathcal{P}(f),f)-
J_2^\tau(\mathcal{P}(g),g)\|_r+\|J_4^\tau(f,\mathcal{P}(f))-
J_4^\tau(g,\mathcal{P}(g))\|_r\right)\\
&\quad\Big.+C_1\left(\|J_1(f,f)-J_1(g,g)\|_r+2\|J_{2}(f,f)-J_{2}(g,g)\|_r
+\|J_4(f,f)-J_4(g,g)\|_r\right)\Big]\\
&\le (1+M\tau)\|f-g\|_r,
\end{align*}
where
\[M=M_2+2M_3+(2+C_1)QD_r.\]
This completes the proof.
\hfill $\square$ \bigskip

Combining the local error bound \eqref{local2} with the stability estimate \eqref{stab2}, and applying a similar argument as in the proof of Theorem \ref{first-th}, we get the error bound of the second-order scheme \eqref{znp} with $u^n$ given by \eqref{sch2} as follows.

\begin{theorem}\label{sec-th}
Let $r>1/2$ and $0<\gamma\le 1$. Assume that the exact solution of \eqref{NF} satisfies $u\in H^{1+r+2\gamma}$ for $0\le t\le T$. Then there exists a constant $\tau_0>0$ such that for all step sizes
$0<\tau\le \tau_0$ and $t_n=n\tau\le T$ we have that the global error of \eqref{sch2} is bounded by
\[\|u(t_n)-u^n\|_r\le C\tau^{1+\gamma},\]
where $C$ depends on $T$, $c$, $r$ and $\|u\|_{L^\infty(0,T;H^{1+r+2\gamma})}$.
\end{theorem}

\begin{corollary}\label{2-th}
Let $r>1/2$ and $0<\gamma\le 1$. Assume that the exact solution of \eqref{GBp} satisfies $z\in H^{1+r+2\gamma}$ and $z_t\in H^{r+2\gamma-1}$ for $0\le t\le T$. Then there exists a constant $\tau_0>0$ such that for all step sizes
$0<\tau\le \tau_0$ and $t_n=n\tau\le T$ we have that the global error of \eqref{znp} combined with \eqref{sch2} is bounded by
\[\|z(t_n)-z^n\|_r+\|z_t(t_n)-z_t^n\|_{r-2}\le C\tau^{1+\gamma},\]
where $C$ depends on $T$, $c$, $r$ and $\|z\|_{L^\infty(0,T;H^{1+r+2\gamma})}+\|z_t\|_{L^\infty(0,T;H^{r+2\gamma-1})}$.
\end{corollary}

Note that the error estimate was established under the constraint $r>1/2$, which enables us to use the the bilinear estimate \eqref{bi} that is crucial for stability (cf. \eqref{stab2}). However, we can derive the error bound in $L^2$ following the approach in \cite{lubich2008, knoller2018fourier} by using a so-called refined bilinear estimate. For fixed $\ep\in (0, 1/2)$, Theorem \ref{sec-th} implies that for solutions in $H^3$, the second-order scheme \eqref{sch2} converges at order $\tau^{3/2}$ in $H^{1}$. This gives an a priori bound on the numerical solution $u^n$ in $H^{1}$. This together with the refined bilinear estimate
\be\label{bi2}
\|fg\|\le \|f\|\|g\|_{L^\infty}\le c\|f\|\|g\|_{1},
\ee
which is a combination of H\"older's inequality and the Sobolev embedding theorem, yields second-order convergence in $L^2$.

\begin{corollary}\label{col2}
Assume that the exact solution of \eqref{GBp} satisfies $z\in H^3$ and $z_t\in H^1$
for $0\le t\le T$. Then there exists a constant $\tau_0>0$ such that for all step sizes
$0<\tau\le \tau_0$ and $t_n=n\tau\le T$ we have that the global error of \eqref{znp} with $u^n$ given by \eqref{sch2} is bounded by
\[\|z(t_n)-z^n\|+\|z_t(t_n)-z_t^n\|_{-2}\le C\tau^2,\]
where $C$ depends on $T$, $c$ and $\|z\|_{L^\infty(0,T;H^3)}+\|z_t\|_{L^\infty(0,T; H^1)}$.
\end{corollary}
\section{Numerical experiments}
In this section, we present some numerical experiments to illustrate our analytic convergence rate given in Corollaries \ref{1-th}, \ref{2-th} and \ref{col2}. In the numerical experiments, we use a standard Fourier pseudospectral method for space discretization where we set the spatial mesh size $\Delta x=1/2^{10}$ in \emph{Example 1} and $\Delta x=\pi/2^{10}$ in \emph{Example 2}, respectively.

\bigskip

{\sl Example 1.} In the first experiment, we study the temporal error of the newly developed exponential-type integrators \eqref{schu} and \eqref{sch2} for the solitary-wave solution. The well-known solitary-wave solution of the GB equation \eqref{GB} is given by
\be\label{sol-ex}
z(x,t)=-\fl{3}{2}\lambda^2\,\sech^2(\lambda/2(x-vt-x_0)),\quad v=\pm (1-\lambda^2)^{1/2},
\ee
where $0<\lambda\le 1$ and $x_0$ are real parameters. Since the solitary wave decays exponentially in the far field, it is possible to use periodic boundary conditions when the domain is chosen large enough. Here we choose $\lambda=1/2$, $x_0=0$ and the torus $\Omega=(-40,40)$.

\begin{figure}[h!]
\begin{minipage}[t]{0.5\linewidth}
\centering
\includegraphics[width=3.15in,height=2.4in]{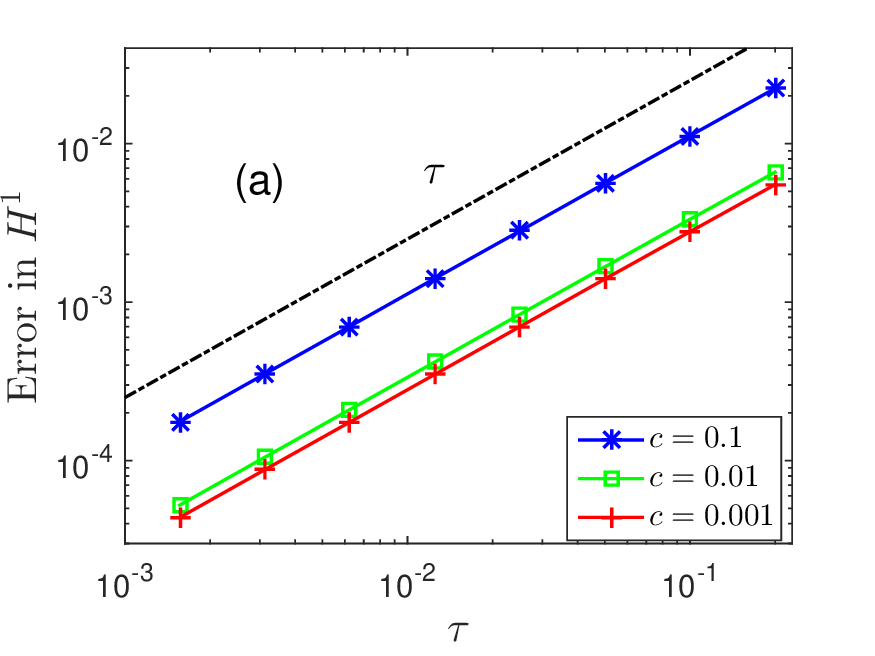}
\end{minipage}%
\hspace{3mm}
\begin{minipage}[t]{0.5\linewidth}
\centering
\includegraphics[width=3.15in,height=2.4in]{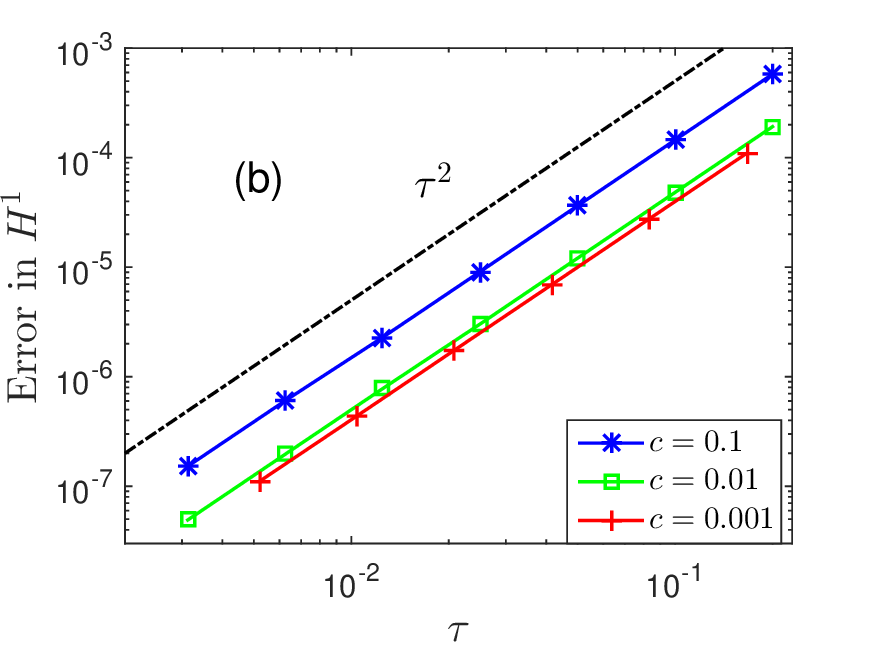}
\end{minipage}
\caption{Numerical simulation for the solitary-wave solution \eqref{sol-ex} at $T=2$ of the first- and second-order exponential-type integrator schemes with spatial mesh size $\Delta x=1/2^{10}$. (a) Linear convergence of the first-order scheme \eqref{schu} for various choices of $c$. The broken line has slope one. (b) Quadratic convergence of the second-order scheme \eqref{sch2} for various choices of $c$. The broken line has slope 2.}\label{fig1}
\end{figure}

Figure \ref{fig1} displays the $H^1$-error of the first- and second-order exponential-type schemes \eqref{schu} and \eqref{sch2}, respectively, at $T=2$ for various $c$ and $\tau$. It can be clearly observed that the schemes \eqref{schu} and \eqref{sch2} converge linearly and quadratically in time, respectively. Moreover, the error decreases as $c$ gets smaller.

\bigskip
{\sl Example 2.} In the second experiment, we apply the newly developed exponential-type integrators \eqref{schu} and \eqref{sch2} to the GB equation with rough initial data. We compare the results with the classical first-order Gautschi-type \cite{hairer2002geometric, grimm2006use} and second-order Deuflhard-type \cite{deuflhard1979study, zhao2016error} exponential integrators used directly to the original GB equation \eqref{GB}.

\begin{figure}[t!]
\begin{minipage}[t]{0.5\linewidth}
\centering
\includegraphics[width=3.15in,height=2.6in]{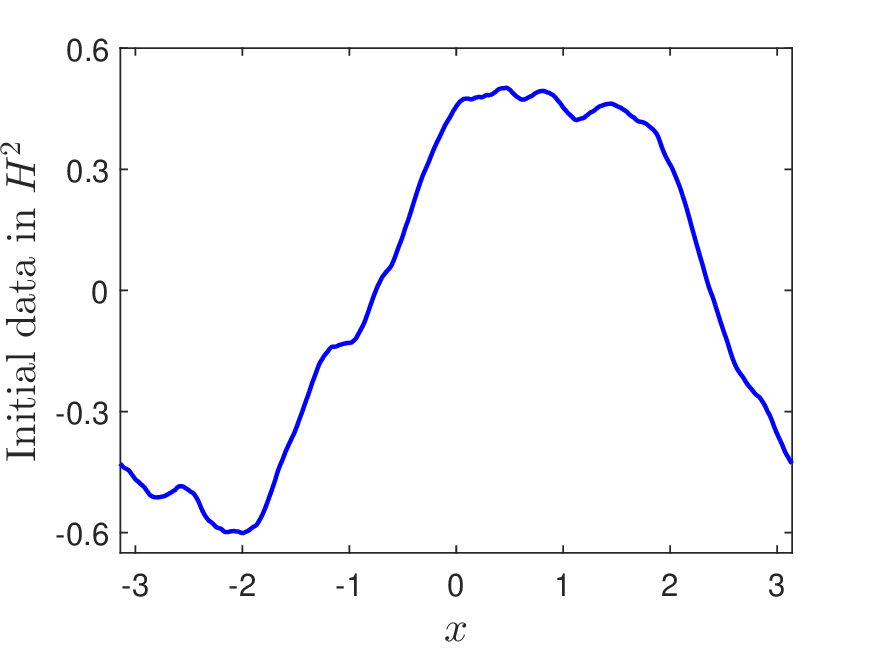}
\end{minipage}%
\hspace{3mm}
\begin{minipage}[t]{0.5\linewidth}
\centering
\includegraphics[width=3.15in,height=2.6in]{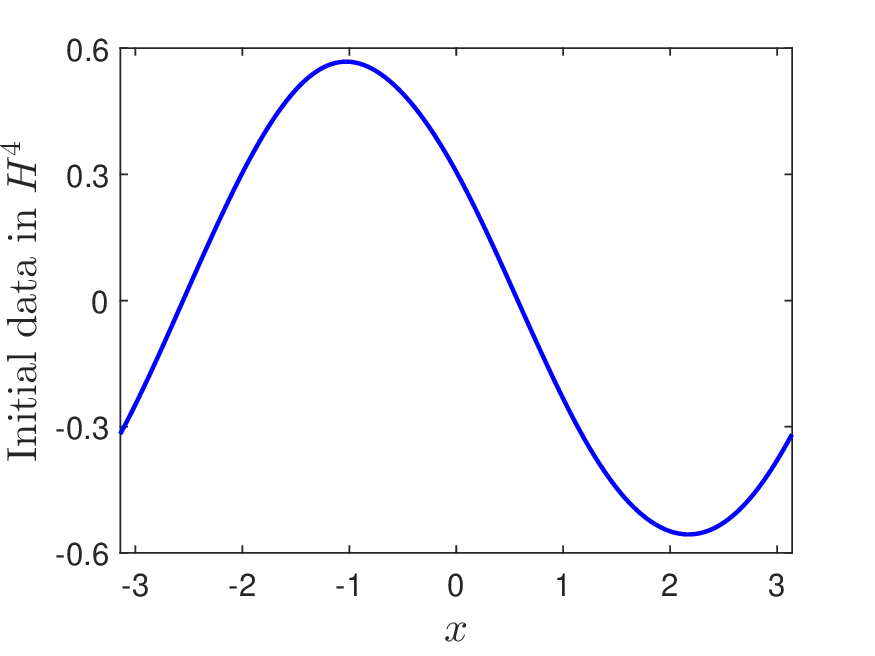}
\end{minipage}
\caption{Initial data $z_0$ normalized in $L^2$ for two different values of $\theta$. Left: $H^2$ initial value for $\theta=2$. Right: $H^4$ initial value for $\theta=4$.}\label{Reg-initial}
\end{figure}

For the grid points $x_j=-\pi+j\Delta x$, $0\le j< M$ for $M=2^{11}$ and $\Delta x=2\pi/M$, we denote by $Z_j$ the approximation of $z(x_j)$. For a vector
\[Z^M:=[Z_0, \cdots, Z_{M-1}]=\mathrm{rand}(1, M)\in\mathbb{R}^M,\]
we set
\[z_\theta^M=|\p_{x,M}|^{-\theta}Z^M,\quad (|\p_{x,M}|^{-\theta})_k:=\left\{\begin{aligned}
&|k|^{-\theta}\quad &\mathrm{if}\quad k\neq 0,\\
&0\quad &\mathrm{if}\quad k=0,
\end{aligned}\right.\]
for different values of $\theta$ normalized in $L^2$ to represent the initial condition in $H^\theta$. For initial conditions with different regularity, we refer to Figure \ref{Reg-initial}.

\begin{figure}[t!]
\begin{minipage}[t]{0.5\linewidth}
\centering
\includegraphics[width=3.15in,height=2.6in]{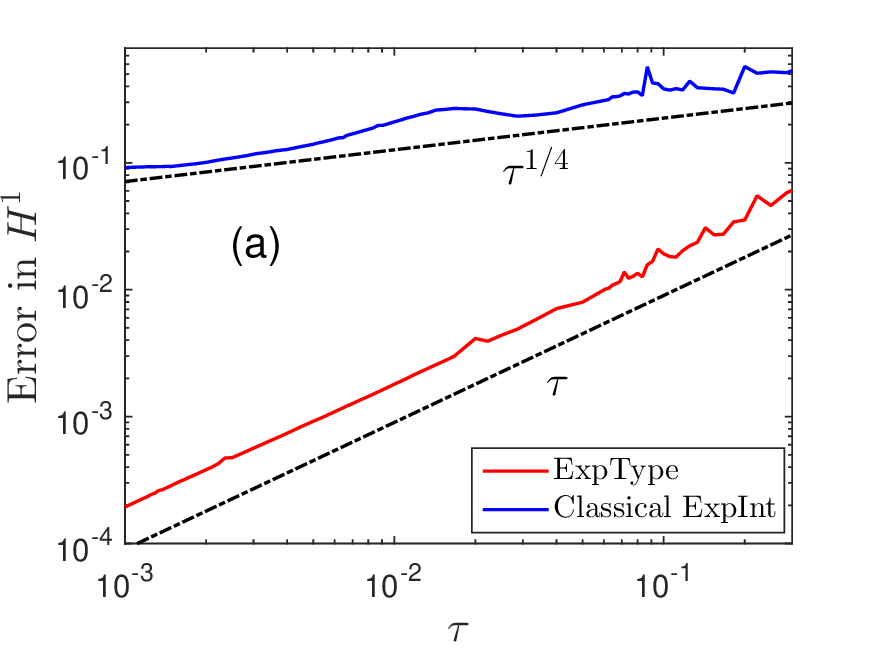}
\end{minipage}%
\hspace{3mm}
\begin{minipage}[t]{0.5\linewidth}
\includegraphics[width=3.15in,height=2.6in]{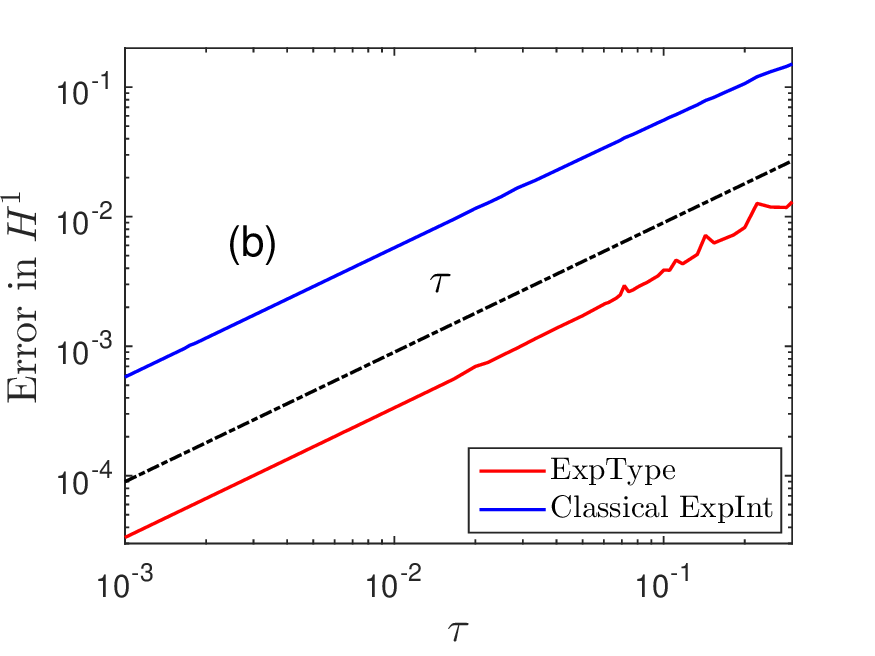}
\end{minipage}
\caption{Temporal error of the first-order exponential-type integrator \eqref{schu} (red) and classical first-order exponential integrator (blue) at $T=2$. The error is measured in $H^1$ with initial data in $H^2$ (a) and $H^3$ (b), respectively. The broken lines are of slope a quarter and one, respectively.}\label{1st-eor}
\end{figure}

In the following, we study the temporal error at time $T=2$ measured in the $L^2$ norm or $H^1$ norm for non-smooth solutions. The parameter $c$ for the newly exponential-type
integrators is chosen as $c=0.01$. The reference solution is obtained by the second-order Deuflhard-type exponential integrator with $\tau=10^{-6}$ and $\Delta x=\pi/2^{10}$.

Figure \ref{1st-eor} displays the convergence behavior of the first-order schemes. The error is measured in $H^1$ with initial data in $H^2$ and $H^3$, respectively. For initial data in $H^2$, the order of the classical first-order Gautschi-type exponential integrator drops to one quarter whereas the newly developed exponential-type integrator is still first-order convergent. The convergence agrees with the error
estimate given in Corollary \ref{1-th}. The first-order convergence of the classical exponential integrator is recovered for smoother initial data in $H^3$. Furthermore, it can be clearly seen from Figure \ref{1st-eor} (b) that the newly exponential type integrator \eqref{schu} is much more accurate than the classical one.

\begin{figure}[h!]
\begin{minipage}[t]{0.5\linewidth}
\centering
\includegraphics[width=3.15in,height=2.6in]{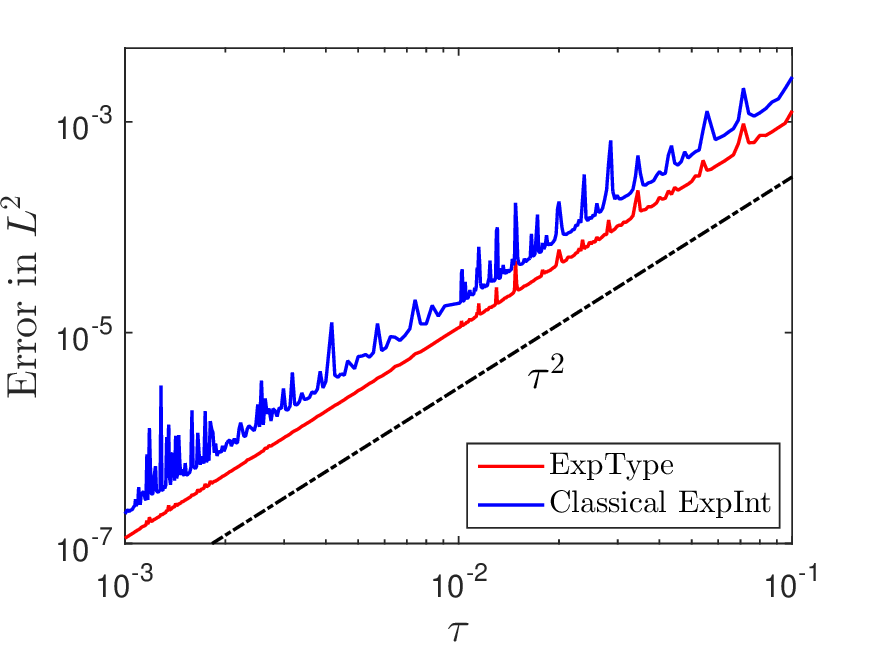}
\end{minipage}%
\hspace{3mm}
\begin{minipage}[t]{0.5\linewidth}
\centering
\includegraphics[width=3.15in,height=2.6in]{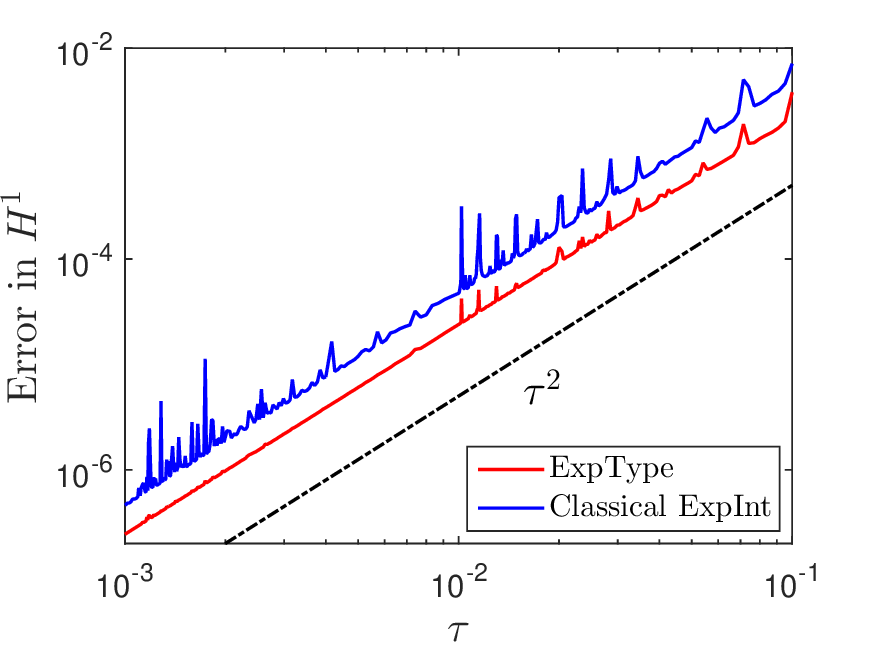}
\end{minipage}
\caption{Temporal error of the second-order exponential-type integrator \eqref{sch2} (red) and classical second-order Deuflhard-type exponential integrator (blue) at $T=2$. Left: the error is measured in $L^2$ with initial data in $H^3$. Right: the error is measured in $H^1$ with initial data in $H^4$. The broken line is of slope two.}\label{2nd-eor}
\end{figure}

\begin{figure}[h!]
\begin{minipage}[t]{0.5\linewidth}
\centering
\includegraphics[width=3.15in,height=2.6in]{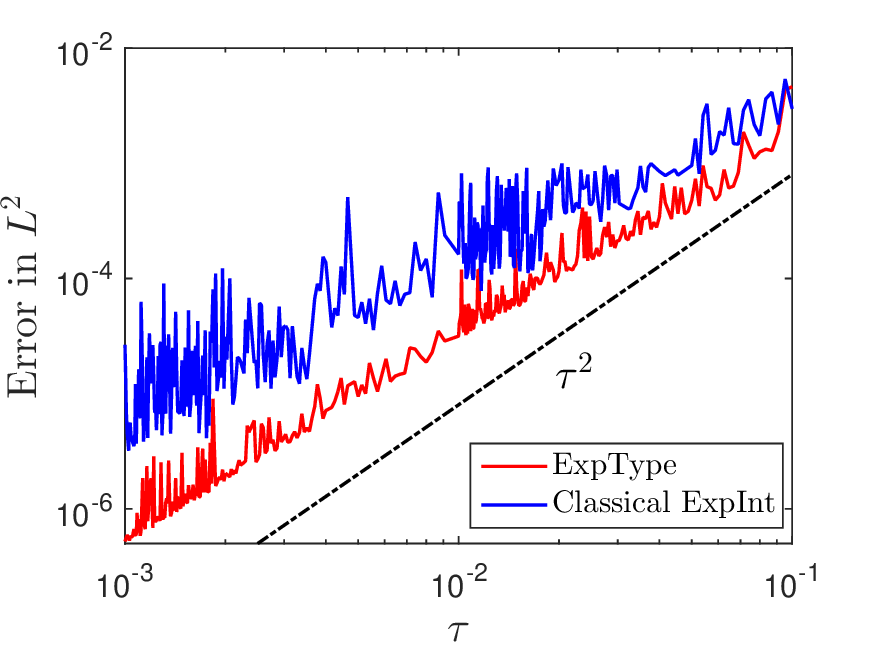}
\end{minipage}%
\hspace{3mm}
\begin{minipage}[t]{0.5\linewidth}
\centering
\includegraphics[width=3.15in,height=2.6in]{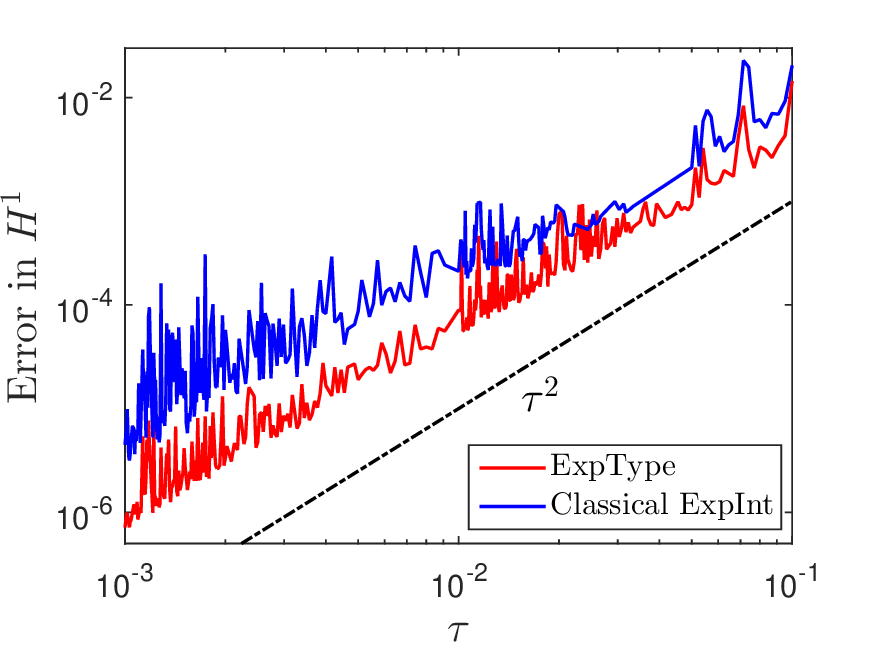}
\end{minipage}
\caption{Temporal error of the second-order exponential-type integrator \eqref{sch2} (red) and classical second-order Deuflhard-type exponential integrator (blue) at $T=2$. Left: the error is measured in $L^2$ with initial data in $H^2$. Right: the error is measured in $H^1$ with initial data in $H^3$. The broken line is of slope two.}\label{ill-sec}
\end{figure}

\begin{figure}[h!]
\begin{minipage}[t]{0.5\linewidth}
\centering
\includegraphics[width=3.15in,height=2.6in]{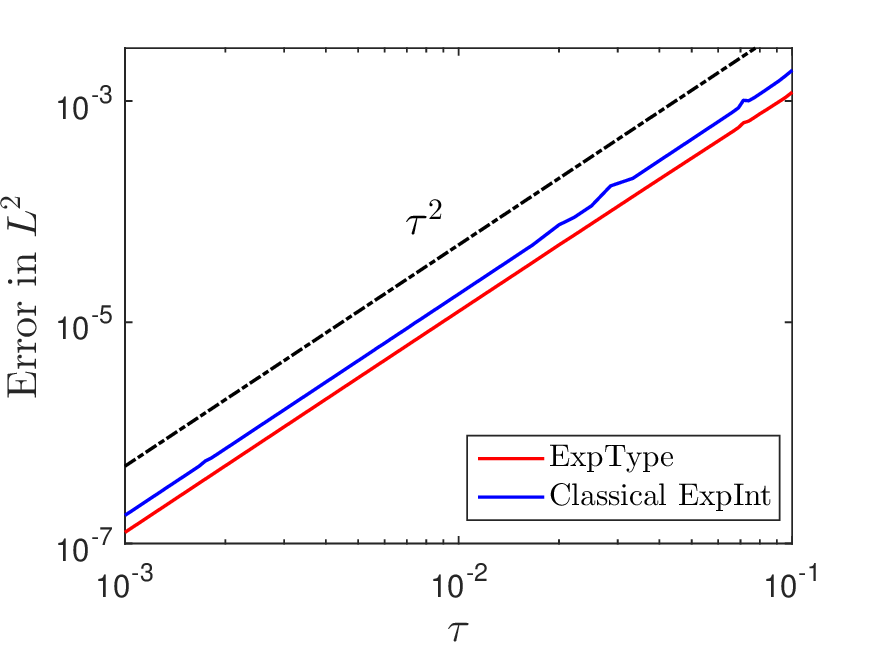}
\end{minipage}%
\hspace{3mm}
\begin{minipage}[t]{0.5\linewidth}
\centering
\includegraphics[width=3.15in,height=2.6in]{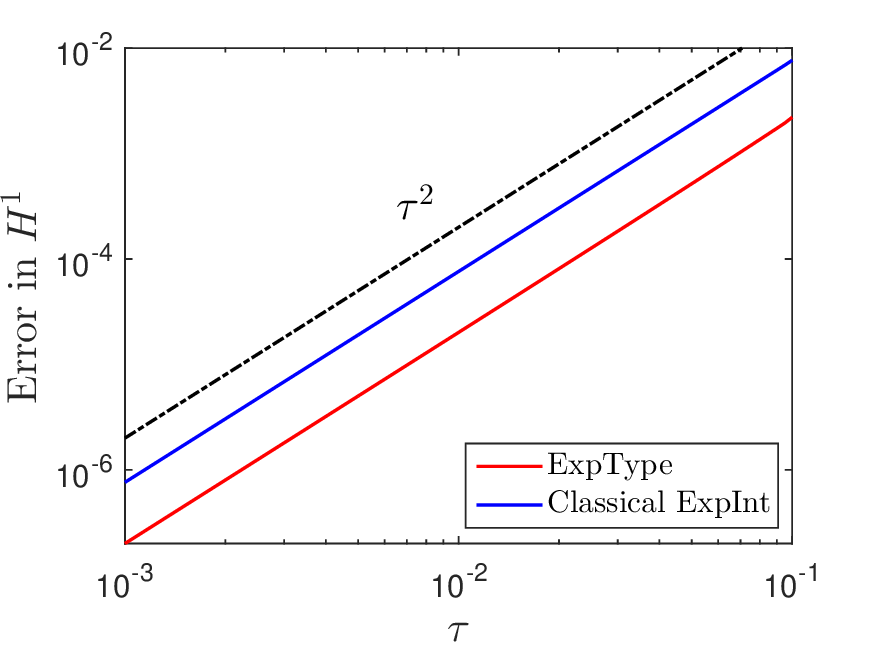}
\end{minipage}
\caption{Temporal error of the second-order exponential-type integrator \eqref{sch2} (red) and classical second-order Deuflhard-type exponential integrator (blue) at $T=2$. Left: the error is measured in $L^2$ with initial data in $H^4$. Right: the error is measured in $H^1$ with initial data in $H^5$. The broken line is of slope two.}\label{Rec-2nd}
\end{figure}

The temporal errors of the second-order schemes are presented in Figure \ref{2nd-eor}: errors in $L^2$ for initial data in $H^3$ (left) and errors in $H^1$ with initial data in $H^4$ (right). For such non-smooth initial conditions, it clearly shows that the
convergence behavior of the second-order Deuflhard-type exponential integrator becomes
irregular while the newly developed second-order exponential-type integrator \eqref{sch2} remains quadratically convergent. This confirms the error estimates in Corollaries \ref{2-th} and \ref{col2}. Figure \ref{ill-sec} shows the errors with less smooth solutions, where the convergence for both exponential integrators behaves irregularly.
This indicates to a certain degree that our error estimates in Corollaries \ref{2-th} and \ref{col2} are sharp. For smoother initial conditions which enable the classical exponential integrator to recover the quadratic convergence rate, we refer to Figure \ref{Rec-2nd}  for error in $L^2$ (left) and error in $H^1$ (right), respectively. Combining the results shown in Figure \ref{1st-eor} (b), we can clearly observe that it requires one additional regularity on the solutions for the classical exponential integrators to get the same convergence rates.

\section{Conclusions}
Two exponential-type integrators were proposed and analyzed for the ``good" Boussinesq (GB) equation with rough initial data. A prior error estimates were established. For order of convergence, it requires one and three additional derivatives for the first- and the second-order exponential-type integrator, respectively. Numerical experiments confirm our analytical results. On the other hand, for non-smooth solutions, the numerical examples show the reliability and superiority of the newly developed exponential-type integrators compared to the classical exponential integrators, which show irregularities and order reductions.

\end{document}